\documentclass{article}

\newtheorem{theorem}{Theorem}
\newtheorem{lemma}{Lemma}

\title{Basics of the theory of cyclic rook polynomials and cyclic
permanents of rectangular matrices}
\author{A.~M.~Kamenetskii}
\date{}

\begin{document}

\maketitle

Let~$A$ be a $k\times m$ matrix, $n\le k$, and let ${j}^{({i})},
{j}^{(i)}_1,\dots,{j}^{({i})}_{p_{{i}}}\in N_m$, $1\le{i}\le q$.
Then, by definition,
$A[N_n|{j}^{({i})}\oplus\sum_{l=1}^{p_{{i}}}\oplus
{j}_l^{({i})}]=A[N_n|{j}^{({i})}]+\sum_{l=1}^{p_{{i}}}A[N_n|
{j}_l^{({i})}]$ and $A[N_n|
{j}^{(1)}\oplus\sum_{l=1}^{p_{1}}\oplus{j}_l^{(1)},\dots,
{j}^{(q)}\oplus\sum_{l=1}^{p_q}\oplus {j}_l^{(q)}]$ is a $n\times
q$ matrix, whose ${i}$~th column is $A[N_n|
{j}^{({i})}\oplus\sum_{l=1}^{p_{{i}}}\oplus{j}_l^{({i})}]$.

\begin{theorem}
Let $0\le r\le n$, $n\ge1$, $k\ge1$, $t\ge0$. Then
\begin{eqnarray*}
\left(\sum_{{i}=0}^ta_{{i}}P_n^{-r+{i}}\right)\otimes J_k&=&
\left(\left( \sum_{{i}=0}^ta_tT_{n+t}^{({i})}\right)\otimes
J_k\right) [N_{nk}|(k(r+1)\oplus\sum_{{1\le{i}\le t\atop
{i}\equiv(r+1){\rm mod}
n}}\oplus k(n+{i}))^{\langle k\rangle},\\
&&(k(r+2)\oplus\sum_{{1\le{i}\le t\atop {i}\equiv(r+2){\rm mod}
n}}\oplus k(n+{i}))^{\langle k\rangle},\dots,
(kn\oplus\sum_{{1\le{i}\le t\atop {i}\equiv n{\rm mod} n}}\oplus k(n+{i}))^{\langle k\rangle},\\
&&(k\oplus\sum_{{1\le{i}\le t\atop {i}\equiv 1{\rm mod} n}}\oplus
k(n+{i}))^{\langle k\rangle}, (2k\oplus\sum_{{1\le{i}\le t\atop
{i}\equiv 2{\rm mod} n}}\oplus k(n+{i}))^{\langle
k\rangle},\dots,\\
&&(kr\oplus\sum_{{1\le{i}\le t\atop {i}\equiv r{\rm mod} n}}\oplus
k(n+{i}))^{\langle k\rangle}].
\end{eqnarray*}
\end{theorem}

An injective mapping~$\varphi$ of a set~$A$ induces the following
bijective mapping~$\varphi^*$ of the set $\varphi(A)\setminus
(\varphi(A)\cap A)$ onto $A\setminus(A\cap\varphi(A))$: if
${j}\in\varphi(A)\setminus(\varphi(A)\cap A)$, then
$\varphi^*({j})={j}_m$, where ${j}_m$ is uniquely determined by
the conditions ${j}_0={j}$ and ${j}_1, {j}_2, \dots,
{j}_{m-1}\in\varphi(A)$, ${j}_l=\varphi^{-1}({j}_{l-1})$, $1\le
l\le m$, ${j}_m\in A\setminus(\varphi(A)\cap A)$. Let
$\bar\alpha=(x_1,\dots,x_n)$ be a sequence, $A\subseteq N_n$,
$\varphi\in {\rm Inj}(A,N_n)$. Denote by $\bar\varphi(\bar\alpha)$
the sequence obtained from  $\bar\alpha$ in the following way. The
components $x_{{j}}$, ${j}\in\varphi(A)$, are deleted. And if
${j}\in\varphi(A)\setminus(\varphi(A)\cap A)$, then the deleted
component~$x_{{j}}$ is replaced by $x_{\varphi^*({j})}$.

\begin{theorem} {\rm\bf (Expansion of the cyclic rook
polynomial of a rectangular matrix along the last~$k$ rows)} Let
$A=(a_{{i},{j}})_{{1\le{i}\le m\atop1\le{j}\le n}}$ be an $m\times
n$ matrix over a commutative ring with unity, $m\le n$, $1\le k\le
m-1$. Then
\begin{eqnarray*}
R(x;z;A)&=&\sum_{S\subseteq N_m\setminus N_{m-k},S\ne\emptyset}
\sum_{\varphi\in {\rm Inj}(S,N_n)} z^{|\varphi|}(\prod_{{i}\in
S}a_{{i},\varphi({i})})R(x;z;A[N_{m-k}|\bar\varphi(N_n)]x^{|S|}\\
&&+R(x;z;A[N_{m-k}|N_n]).
\end{eqnarray*}
\end{theorem}

If $k>1$, then such an expansion along arbitrary~$k$ rows, for an
arbitrary matrix, does not exist. On the other hand, if $k=1$,
then such an expansion along arbitrary row exists.

\begin{theorem}
Let $A=(a_{{i},{j}})_{{1\le{i}\le m\atop1\le{j}\le n}}$ and $m\le
n$, $1\le i\le m$. Then
\begin{eqnarray*}
R(x;z;A)&=&a_{{i},{i}}xzR(x;z;A[N_m\setminus({i})|N_n\setminus({i})]\\
&&+\sum_{{j}\in
N_n\setminus({i})}a_{{i},{j}}xR(x;z;A[N_m\setminus({i})|1,2,\dots,{i}-1,
{i}+1,\dots,{j}-1,{i},{j}+1,\dots,n]\\&&+
R(x;z;A[N_m\setminus({i})|N_m\setminus({i}),{i},m+1,\dots,n].
\end{eqnarray*}
\end{theorem}

\begin{theorem}{\rm\bf (Expansion of the cyclic rook
polynomial along arbitrary~$k$ rows)} Let
$A=(a_{{i},{j}})_{{1\le{i}\le m\atop1\le{j}\le n}}$ ¨ $m\le n$,
$1\le k\le m-1$, $\bar\beta\in Q_{k,m}$. Then
$$
{\rm per}(z;A)=\sum_{\varphi\in{\rm Inj}(\{\bar\beta\},N_n)}
z^{|\varphi|}(\prod_{i\in\{\bar\beta\}}a_{i,\varphi(i)}){\rm
per}(z;A[N_m\setminus\bar\beta|\bar\varphi(N_n)].
$$
\end{theorem}

\begin{theorem}
Let $A=(a_{i,j})$ and $B=(b_{i,j})$ be $m\times n$ matrices over a
commutative ring with unity, $m\le n$,
$R(x;z;A)=\sum_{l=0}^mr_l(z;A)x^l$, $r_0(z;A)=1$,
 $R(x;A)=\sum_{l=0}^mr_l(A)x^l$, $r_0(A)=1$. Then
\begin{eqnarray*}
R(x;z;A+B)&=&\sum_{s=0}^m\sum_{\bar\alpha\in Q_{s,m}}
\sum_{\varphi\in {\rm Inj}(\{\bar\alpha\}, N_n)}
z^{|\varphi|}\left(\prod_{i\in\{\bar\alpha\}}a_{i,\varphi(i)}\right)
x^sR(x;z;B[N_m\setminus\bar\alpha|\bar\varphi(N_n)]);\\
r_l(z;A+B)&=&\sum_{s=0}^l\sum_{\bar\alpha\in Q_{s,m}}
\sum_{\varphi\in {\rm Inj}(\{\bar\alpha\}, N_n)}
z^{|\varphi|}\left(\prod_{i\in\{\bar\alpha\}}a_{i,\varphi(i)}\right)
r_{l-s}(z;B[N_m\setminus\bar\alpha|\bar\varphi(N_n)]);\\
{\rm per}(z;A+B)&=&\sum_{s=0}^m\sum_{\bar\alpha\in Q_{s,m}}
\sum_{\varphi\in {\rm Inj}(\{\bar\alpha\}, N_n)}
z^{|\varphi|}\left(\prod_{i\in\{\bar\alpha\}}a_{i,\varphi(i)}\right)
{\rm per}(z;B[N_m\setminus\bar\alpha|\bar\varphi(N_n)]);\\
R(x;A+B)&=&\sum_{s=0}^m\sum_{{\bar\alpha\in Q_{s,m}\atop\bar\beta\in Q_{s,n}}}
{\rm per}(A[\bar\alpha|\bar\beta])x^sR(x;B(\bar\alpha|\bar\beta));\\
r_l(A+B)&=&\sum_{s=0}^l\sum_{{\bar\alpha\in Q_{s,m}\atop\bar\beta\in Q_{s,n}}}
{\rm per}(A[\bar\alpha|\bar\beta])r_{l-s}(B(\bar\alpha|\bar\beta));\\
{\rm per}(A+B)&=&\sum_{s=0}^m\sum_{{\bar\alpha\in Q_{s,m}\atop\bar\beta\in Q_{s,n}}}
{\rm per}(A[\bar\alpha|\bar\beta]){\rm per}(B(\bar\alpha|\bar\beta)).
\end{eqnarray*}
\end{theorem}

\begin{theorem}
Let~$A$ be an $m\times n$ matrix over a commutative ring with
unity, $m\le n$, $R(x;z;A)=\sum_{l=0}^m r_l(z;A)x^l$,
$r_0(z;A)=1$, $(z)^{(k)}=\sum_{{i}=0}^{k-1}(z+{i})$, $k\ge1$,
$(z)^{(0)}=1$. Then
$$
r_l(z;yJ_{m,n}-A)=\sum_{s=0}^l(-1)^s{m-s\choose
l-s}r_s(z;A)(n-l+z)^{(l-s)}y^{l-s},\qquad 0\le l\le m.
$$
In particular,
$$
{\rm
per}(z;yJ_{m,n}-A)=\sum_{s=0}^m(-1)^sr_s(z;A)(n-m+z)^{(m-s)}y^{m-s}.
$$
\end{theorem}

\begin{theorem}
$$
{\rm per}(z;(a_0I_n+a_1P_n)\otimes
J_k)=\sum_{s=0}^k\frac{(z)^{(s)}} {s!}[s!{k\choose
s}a_0^{k-s}a_1^s(s+z)^{(k-s)}]^n
$$
for all $n\ge1$.
\end{theorem}

\begin{theorem}
Let $r\ge t+1$. Then for all $n\ge1$ we have
$$
R(x;z;(\sum_{{i}=0}^ta_{{i}}P_n^{-r+{i}})\otimes J_k)=
\sum_{{(l_1,\dots,l_t)\atop0\le l_i\le k,1\le i\le t}} \sum_M
f(z;{v_1^{\langle l_1\rangle},v_2^{\langle l_2\rangle},\dots,
v_t^{\langle l_t\rangle}, v_{t+1}^{\langle k\rangle},
v_{t+2}^{\langle k\rangle},\dots,v_r^{\langle k\rangle}\choose
\bar\gamma\setminus(1^{\langle k-l_1\rangle},2^{\langle
k-l_2\rangle},\dots,t^{\langle k-l_t\rangle})}\cdot
$$
$$
(K_{[k],t}^{\langle0,r\rangle}(a_0x,a_1x,\dots,a_tx;z))^n
[L(1^{\langle k-l_1\rangle},v_1^{\langle l_1\rangle},2^{\langle
k-l_2\rangle},v_2^{\langle l_2\rangle},\dots,t^{\langle
k-l_t\rangle},v_t^{\langle l_t\rangle},v_{t+1}^{\langle k\rangle},
v_{t+2}^{\langle k\rangle},\dots,v_r^{\langle
k\rangle})|L(\bar\gamma)],
$$
where the summation is carried over the domain~$M$ defined by the
conditions
\begin{eqnarray*}
\bar\gamma&\in&\cup_{d=0}^{rk}G_{[k],d,t}^{\langle0,r\rangle},\\
\{\bar\gamma\}\setminus(\{\bar\gamma\}\cap\{w^{\langle
rk\rangle}\})&=&\{1^{\langle k-l_1\rangle},v_1^{\langle
l_1\rangle},2^{\langle k-l_2\rangle},v_2^{\langle
l_2\rangle},\dots,t^{\langle k-l_t\rangle},v_t^{\langle
l_t\rangle},v_{t+1}^{\langle k\rangle}, v_{t+2}^{\langle
k\rangle},\dots,v_r^{\langle k\rangle}\}.
\end{eqnarray*}
\end{theorem}

\begin{lemma}
Let $\bar\alpha=(x_1,\dots,x_m)$ and $\bar\gamma=(y_1,\dots,y_n)$
be two sequences, $m\le n$,
$\{\bar\alpha\}\subseteq\{\bar\gamma\}$,
$\{\bar\alpha\}=\{b_1^{\langle l_1\rangle},\dots,b_d^{\langle
l_d\rangle}\}$, $A_i=\{j\in N_m|x_j=b_i\}$, $B_i=\{j\in
N_n|y_j=b_i\}$, $1\le i\le d$. Then, for all $i, 1\le i\le d$, we
have
$$
f(z;{\bar\alpha\choose\bar\gamma})=
({|B_i|\choose|A_i|}(|A_i|!))^{-1} \sum_{\varphi\in{\rm
Inj}(A_i,B_i)} z^{|\varphi|}f(z;{\bar\alpha\setminus(b_i^{\langle
l_i\rangle})\choose \bar\varphi(\bar\gamma)}).
$$
\end{lemma}

The author expresses his deep gratitude to S.~K.Lando for
invaluable support and attention to the work.

{\tt e-mail: kamenetsky\_a@list.ru}
\end{document}